\documentclass[12pt,leqno]{article}
\usepackage{amsmath,amsfonts,amssymb,color,hyperref}
\usepackage{color}

\numberwithin{equation}{section}

\newtheorem{theorem}{Theorem}[section]
\newtheorem{definition}[theorem]{Definition}
\newtheorem{proposition}[theorem]{Proposition}
\newtheorem{corollary}[theorem]{Corollary}
\newtheorem{lemma}[theorem]{Lemma}
\newtheorem{remark}[theorem]{Remark}
\newtheorem{example}[theorem]{Example}

\newcommand{\ang}[1]{\langle  #1 \rangle }  

\newcommand{\p}{\mathbb{P}}
\newcommand{\Q}{\mathbb{Q}}

\newcommand{\E}[1]{\mathbb{E}\left[#1\right]}
\newcommand{\Et}[1]{\mathbb{E}_t\left[#1\right]}

\newcommand{\trieq}{\stackrel{\triangle}{=}}
\newcommand{\eof}{\hfill {\it Q.E.D.} \vspace*{0.3cm}}

\newcommand{\pf}{{\it Proof: }}

\newcommand{\cF}{{\mathcal F}}
\newcommand{\cG}{{\mathcal G}}

\newcommand{\R}{\mathbb{R}}
\newcommand{\id}{{\mathbf 1}}

\title{Equilibrium for Time-Inconsistent Stochastic Linear--Quadratic Control under Constraint}

\author{Ying Hu\thanks{IRMAR,
Universit\'e Rennes 1, 35042 Rennes Cedex, France (ying.hu@univ-rennes1.fr)  and School of
Mathematical Sciences, Fudan University, Shanghai 200433, China, partially supported by Lebesgue center of mathematics ``Investissements d'avenir"
program - ANR-11-LABX-0020-01,  by ANR CAESARS - ANR-15-CE05-0024 and by ANR MFG - ANR-16-CE40-0015-01.} \and
Jianhui Huang\thanks{Department of Applied Mathematics, The Hong Kong Polytechnic University, Hong Kong (majhuang@polyu.edu.hk), partially supported by Hong Kong RGC Grant 502412,
15300514, G-YL04. } \and
Xun Li\thanks{Department of Applied Mathematics, The Hong Kong Polytechnic University, Hong Kong (malixun@polyu.edu.hk), partially supported by PolyU G-UA7U, Hong Kong RGC grants 15224215 and 15255416.}}

\begin{document}
\maketitle

\begin{abstract}
In this paper, we study a class of stochastic time-inconsistent linear-quadratic (LQ) control problems with control input constraints.
These problems are investigated within the more general framework associated with random coefficients.
This paper aims to further develop a new methodology, which fundamentally differs from those in the standard control (without constraints) theory in the literature, to cope with the mathematical difficulties raised due to the presence of input constraints. We first prove that the existence of an equilibrium solution is equivalent to the existence of a solution to some forward-backward stochastic differential equations with constraints.
Under convex cone constraint, an explicit solution to equilibrium for mean-variance portfolio selection can be obtained and proved to be unique. Finally, some examples are discussed to shed light on the comparison between our established results and standard control theory. \end{abstract}


{\bf AMS subject.}
 91B51, 93E99, 60H10

\section{Introduction}
In dynamic decision making, the presence of time inconsistency is often identified in socioeconomic systems and accordingly, its study has important values in various fields, such as engineering, management science, finance and economics (for example, see Kydland and Prescott \cite{KP}). More recently, considerable research attention has been paid in studying this family of stochastic time-inconsistent control problems as well as their financial applications. The study on time-inconsistency by economists can be traced back to Strotz \cite{Strotz} in the 1950s, who initiated the formulation of time-inconsistent decision making as a game between incarnations of the decision maker  himself. For the sake of motivation and to make our discussion concrete, let us briefly lay out a simple but illustrating example of time inconsistency in dynamic setting.
\begin{example}\rm
Let $\lambda > 0$ be a constant. Then we consider the following dynamic mean-variance portfolio problem
%
%
\begin{equation}\label{eq:intro}
\left\{\begin{array}{cl}
\displaystyle\min_{u} & J(u):=\mbox{\rm Var}(X_T)-2\lambda\mathbb{E}[X_T], \\
\mbox{s.t.} & dX_s=\theta'u_sds+u_s'dW_s, \quad\quad s \in [0,T], \\
& ~ X_0 = x_0,
\end{array}\right.
\end{equation}
where $X_\cdot \in \mathbb{R}$, $u_\cdot \in \mathbb{R}^2$, $\theta=(1,1)'$ and $W_\cdot$ is a two-dimensional standard Weiner process.
%
Its dynamic counterpart yields the following optimization problem over $[t,T]$, for any $t\in[0,T]$,
\begin{equation}\label{eq:intro_t}
\left\{\begin{array}{cl}
\displaystyle\min_{u} & J_t(u):=\mbox{\rm Var}_t(X_T)-2\lambda\mathbb{E}_t[X_T], \\
\mbox{s.t.} & dX_s=\theta'u_sds+u_s'dW_s, \quad\quad s \in [t,T], \\
& ~ X_t = x_t,
\end{array}\right.
\end{equation}
where $\Et{\cdot}=\E{\cdot|\cF_t}$ is the conditional expectation and $\mbox{\rm Var}_t(\cdot)$ is the conditional variance under $\Et{\cdot}$.

An admissible control $u_{\cdot}^{*,0,x_0}$ is optimal for the problem (\ref{eq:intro}) if $J(u^{*,0,x_0}_{\cdot}) = \displaystyle\min_{u}J(u)$.
Also, we define the optimal control $u^{*,t,x_t}_\cdot$ for the problem (\ref{eq:intro_t}) similarly.
We say the problem (\ref{eq:intro_t}) is time-consistent if, for any $t \in [0, T]$, it holds that
\begin{equation}\label{eq:intro_u}
u_s^{*,t,X_t^{*,0,x_0}}=u_s^{*,0,x_0} \quad\mbox{for } t \leq s \leq T.
\end{equation}
However, applying results obtained in \cite{HJZ,HZ,LX,LZL}, we have an optimal control
$u_s^{*,0,x_0}=(x_0-X_s^{*,0,x_0}+\lambda e^{2T})\theta$, $0 \leq s \leq T$ for Problem (\ref{eq:intro}), where
$$
X_s^{*,0,x_0}=x_0+2\int_0^s(x_0-X_v^{*,0,x_0}+\lambda e^{2T})dv+\int_0^s(x_0-X_v^{*,0,x_0}+\lambda e^{2T})\theta'dW_s.
$$
Again, as before we should have an optimal control on $[t,T]$, $u_s^{*,t,X_t^{*,0,x_0}}=(X_t^{*,0,x_0}-X_s^{*,t,X_t^{*,0,x_0}}+\lambda e^{2(T-t)})\theta$, $t \leq s \leq T$ for Problem (\ref{eq:intro_t}), where
$$
\begin{array}{rl}
X_s^{*,t,X_t^{*,0,x_0}}=\!\!\! & X_t^{*,0,x_0}+\displaystyle2\int_t^s(X_t^{*,0,x_0}-X_v^{*,t,X_t^{*,0,x_0}}+\lambda e^{2(T-t)})dv \\
& +\displaystyle\int_t^s(X_t^{*,0,x_0}-X_v^{*,t,X_t^{*,0,x_0}}+\lambda e^{2(T-t)})\theta'dW_v. \\
\end{array}
$$
It is obvious that $u_s^{*,t,X_t^{*,0,x_0}} \neq u_s^{*,0,x_0}$ for $ t \leq s \leq T$.
The dynamic optimization problem (\ref{eq:intro_t}) is called time-inconsistent since (\ref{eq:intro_u}) fails to hold. Therefore, time inconsistency reflects that an optimal strategy at present may no longer be optimal in the future.
\hfill $\Box$
\end{example}

In response, Strotz suggested two possible fundamental schemes to circumvent time inconsistency: (i) ``He may try to precommit his future activities either irrevocably or  by contriving a penalty for his future self if he should misbehave", which is named the {\it strategy of pre-commitment}; and (ii) ``He may resign himself to the fact of intertemporal conflict and decide that his `optimal' plan at any date is a will-o'-the-wisp which cannot be attained, and learn to select the present action which will be best in the light of future disobedience", which is termed the {\it strategy of consistent planning}. The strategy of consistent planning is also called the {\it time-consistent policy} in the literature. For a dynamic mean-variance model, Basak and Chabakauri \cite{Basak} reformulated it as an intrapersonal game model where the investor optimally elicits the policy at any time $t$, on the premise that he has already decided time-consistent (equilibrium) policies applied in the future.

The game formulation is tractable to capture time inconsistency when the underlying time setting is (finite or countable) discrete. Nevertheless, when the time setting is continuous, the formulation should be generalized or modified in different ways. Additionally, some tailor-made arguments, to be shown later, should also be introduced to handle the continuous-time setting. We remark that it is still unclear which is the best one among different definitions of a solution to time-inconsistent decision problem. Mathematically,
both the existence and the uniqueness of a solution make a definition more acceptable. Although it is common that a game problem admits multiple solutions, the time-inconsistent decision problem is a decision problem for single player, and hence, an identical value process for all solutions
is considered to be more reasonable even if the control may allow multiple solutions. Instead of seeking an ``optimal control", some kind of equilibrium controls are worthy to be developed in both theoretical methodology and numerical computation algorithm. This is mainly motivated by practical applications in statistical economics and has recently attracted considerable interest and attempts.

Yong \cite{Yong} and Ekeland and Pirvu \cite{EP} established the existence of equilibrium solutions, with their own  definitions for equilibrium solutions, for the time inconsistency caused by hyperbolic discounting. Grenadier and Wang \cite{GW} also studied the hyperbolic discounting problem in an optimal stopping model.
In a Markovian system, Bj\"ork and Murgoci \cite{BM} proposed a definition of a general stochastic control problem with time-inconsistent terms, and
presented some sufficient condition for a control to be a solution by a system of partial differential equations. They constructed some solutions for some examples including an LQ  one, but it looks very hard to find not-too-harsh
condition on parameters to ensure the existence of a solution. Bj\"ork, Murgoci and Zhou \cite{BMZ} also derived an equilibrium for a mean-variance portfolio selection with state-dependent risk aversion.
Basak and Chabakauri \cite{Basak} studied an equilibrium strategy for a mean-variance portfolio selection problem with constant risk aversion and got more details
on the constructed solution.
Hu, Jin and Zhou in \cite{HJZ} generalized the discrete-time game formulation for an LQ  control problem with time-inconsistent terms in a non-Markovian system, which is slightly different from the one in Bj\"ork and Murgoci \cite{BM}, and constructed an equilibrium strategy for quite general LQ control problem, including a non-Markovian system, and then in \cite{HJZ2}, they proved that the constructed equilibrium strategy  is unique. Bensoussan, Frehse and Yam \cite{BFY} introduced a class of time-inconsistent game problems of mean-field type and provided their equilibrium solutions; Karnam, Ma and Zhang \cite{KMZ} introduced the idea of ``dynamic utility" under which the original time-inconsistent problem (under the fixed utility) is transferred to time-consistent one.
In addition, Cui, Li, Wang and Zhu \cite{CLWZ} showed that the multi-period mean-variance
problem does not satisfy time consistency in efficiency and developed a revised mean-variance strategy.
By relaxing the self-financing restriction to allow the withdrawal of money from the market, the revised mean-variance strategy dominates the original dynamic mean-variance strategy in the mean-variance space.
Furthermore, Cui, Li, Li and Shi \cite{CLLS} further investigated the time-consistent strategy for a behavioral
risk aversion model by solving a nested mean-variance game formulation.

Recently, Bensoussan, Wong, Yam and Yung \cite{BWYY} studied the time-consistent strategies in the mean-variance portfolio selection with short-selling prohibition in both discrete-time and continuous-time settings and showed that the discrete-time equilibrium controls converge to that in the continuous-time setting. In their work, the cost functional just includes the terminal mean and variance terms without the running cost part.
In this paper, we further consider a class of time-inconsistent stochastic LQ control problems under control constraint involving the integral part in the cost functional. Also, we investigate these problems within the framework of random coefficients.
Hu, Jin and Zhou \cite{HJZ,HJZ2} introduced the new methodology which distinguishes significantly from those in classic control (without constraints) theory in the literature, to tackle time-inconsistent stochastic LQ control problem without constraints.
Our work aims to further develop the new methodology proposed in \cite{HJZ} to cope with the mathematical difficulties rooted in the presence of control constraints.
We first prove that the existence of an equilibrium solution is equivalent to the existence of a solution to some forward-backward stochastic differential equations (FBSDE) with constraints.
Then we present an explicit solution to equilibrium for mean-variance portfolio selection under convex cone constraint and show that the  constructed solution is unique. Finally, we illustrate the established results using examples.
In particular, we compare our results with that in Bensoussan, Wong, Yam and Yung \cite{BWYY} for the deterministic coefficients. Our current work is one further step toward understanding the role of input constraint in time-inconsistency decision making, and we expect to see more research progress along this direction.

The rest of this paper is organized as follows. In Section \ref{prob-formu}, we give the formulation of the LQ control problem without time consistency under constraint.
Then we give an equivalent characterization of a solution by a system of forward-backward stochastic differential equations in Section \ref{formal-derivation}. Finally in Section \ref{mvcase}, we give an explicit solution to equilibrium for mean-variance portfolio selection under convex cone constraint and show that the thus constructed solution is unique.

\section{Problem Formulation}\label{prob-formu}

Let $(W_t)_{0 \le t \le T}=(W_t^1,\cdots,W_t^d)_{0 \le t \le T}$
be a $d$-dimensional Brownian motion on a probability
space $(\Omega, \cF, \p)$. Denote by $(\cF_t)$ the
augmented filtration generated by $(W_t)$.

We will use the following notation. Let $p\ge 1$.
\begin{equation*}
\begin{tabular}{rl}
$\mathbb S^l$: & the set of symmetric $l \times l$ real matrices. \\
 $L^p_{\cG}(\Omega; \, \R^l)$: & the set of random variables $\xi: (\Omega, \cG) \rightarrow (\R^l, {\cal B}(\R^l))$ with $\E{|\xi|^p}<+\infty$. \\
 $L^\infty_{\cG}(\Omega; \, \R^l)$: & the set of essentially bounded  random variables $\xi: (\Omega, {\cG}) \rightarrow (\R^l, {\cal B}(\R^l))$. \\
 $L^p_\cG(t, \, T; \, \R^l)$: & the set of $\{\cG_s\}_{s\in [t,T]}$-adapted processes $f=\{f_s: t\leq s\leq T\}$ with \\
 & $\E{ \int_t^T|f_s|^p\, ds} < \infty$. \\
 $L^\infty_\cG(t, \, T; \, \R^l)$: & the set of essentially  bounded $\{\cG_s\}_{s\in [t,T]}$-adapted processes. \\
 $L^p_\cG(\Omega; \, C(t, \, T; \, \R^l))$: & the set of continuous $\{\cG_s\}_{s\in [t,T]}$-adapted processes \\
         & $f=\{f_s: t\leq s\leq T\}$ with $\E{ \sup_{s\in [t,T]}|f_s|^p\, } < \infty$.
\end{tabular}
\end{equation*}
%
%
%
%
%
%

We will often use vectors and matrices in this paper, where all vectors are column vectors. For a matrix $M$,
 define
$M'$ as transpose and $\vert M\vert=\sqrt{\sum_{i,j}m_{ij}^2}$ as Frobenius norm of a matrix $M$, respectively.

Now we introduce the model under consideration in this paper.

Let $T>0$ be given and fixed. The controlled system is governed by the following stochastic
differential equation (SDE)
on $[0, T]$:
\begin{equation}\label{controlgeneral:0}
dX_s=[A_sX_s+B_s'u_s+b_s]ds+ \sum_{j=1}^d[C_s^jX_s+D_s^{j} u_s+\sigma_s^{j} ]dW_s^{j}, \quad X_0=x_0,
\end{equation}
where  $A$ is a bounded deterministic  function on $[0, T]$ with values in $\R^{n\times n}$,
$B,C^j,D^j$ are all essentially bounded adapted processes on $[0,T]$
with values in $\R^{ l\times n}$,
$\R^{n\times n}$, $\R^{ n\times l}$,  respectively, and  $b$ and $\sigma^j$ are stochastic processes
in  $L^\infty_\cF(0,T; \R^n)$. Let $K$ be a given convex set in $\mathbb R^l$.
The process $u \displaystyle\in \bigcup_{p> 2} L^p_\cF(0, \, T; \, K)$ is the control, and
$X\in L^p_\cF(\Omega; \, C(0, \, T; \, \R^n))$ is the corresponding state process  with initial value $x_0\in \R^n$ and with $u\in L^p_\cF(0, \, T; \, K)$.

When time evolves to $t\in[0,T]$, we need to  consider the controlled system
starting from  $t$ and state $x_t\in L^p_{\cF_t}(\Omega; \, \R^n)$:
\begin{equation}\label{controlgeneral:t}
dX_s=[A_sX_s+B_s'u_s+b_s]ds+ \sum_{j=1}^d[C_s^jX_s+D_s^{j} u_s+\sigma_s^{j} ]dW_s^{j}, \quad X_t=x_t.
\end{equation}
For any  control $u\in L^p_\cF(t, T;  K)$, there exists
a unique  solution $X^{t,x_t,u}\in L^p_\cF(\Omega; \, C(t,  T;  \R^n))$.
At time $t$ with the system state $X_t=x_t$, our aim  is to minimize
\begin{eqnarray}\label{costgeneral}
J(t,x_t;u)&\trieq&\frac{1}{2}\mathbb E_t\int_t^T\left[ \ang{Q_sX_s, X_s}+\ang{ R_su_s, u_s}\right]ds+\frac{1}{2}\mathbb E_t [\ang{G X_T, X_T}]\nonumber \\
&&- \frac{1}{2} \ang{h\Et{X_T}, \Et{X_T}}-\ang{ \mu_1 x_t+\mu_2,  \Et{X_T}}
\end{eqnarray}
over  $u\in L^p_\cF(t, \, T; \, K)$,
where $X=X^{t,x_t,u}$, and $\Et{\cdot}=\E{\cdot|\cF_t}$.
In the above $Q$ and $R$ are both positive semi-definite and essentially bounded adapted processes on $[0,T]$ with values in ${\mathbb S}^n$
and ${\mathbb S}^l$ respectively,  $G, h, \mu_1, \mu_2$ are constants in $\mathbb S^n$, $\mathbb S^n$,   $\R^{n\times n}$, $\R^n$ respectively, and moreover $G$ is  positive semi-definite.

We define an {\it equilibrium (control)} in the following manner.
Given a  control $u^*$, for any $t\in [0,T)$, $\varepsilon>0$ and $v\displaystyle\in \bigcup_{p>2} L^p_\cF(t, \, T; \, K)$,  define
\begin{equation}\label{svgeneral}
u^{t,\varepsilon,v}_s=u^*_s+ (v_s-u_s^*)\id_{s\in [t,t+\varepsilon)},\;\;\;s\in[t,T].
\end{equation}

\begin{definition}
Let $u^* \in \displaystyle\bigcup_{p>2} L^p_\cF(0, \, T; \, K)$ be a given control and $X^*$ be the state process corresponding to $u^*$. The control $u^*$ is called an equilibrium if
$$
\liminf_{\varepsilon\downarrow 0}\frac{J(t,X^*_t;u^{t,\varepsilon,v})-J(t,X^*_t;u^*)}{\varepsilon}\ge 0,
$$
where $u^{t,\varepsilon,v}$ is defined by (\ref{svgeneral}), for any $t\in [0,T)$ and $v \in \displaystyle\bigcup_{p>2} L^p_\cF(t, \, T; \, K)$.
\end{definition}

\begin{remark} \rm
There is some difference between our definition and that of \cite{HJZ},
because there is a control constraint $K$ in our situation. Note that the convexity of $K$ is not needed in our definition.
\end{remark}

\section{Necessary and Sufficient Condition  of Equilibrium Controls}\label{formal-derivation}

 In this section, we present a general necessary and sufficient condition for equilibria. This condition is made possible by a stochastic Lebesgue differentiation theorem involving conditional expectation.

To proceed, we start with some relevant known result from  \cite{HJZ}.
Let $u^*$ be a fixed control and $X^*$ be the corresponding state process.
For any $t\in [0, T)$, define in the time interval $[t, T]$  the  processes
$(p(\cdot;t), (k^j(\cdot;t))_{j=1,\cdots, d})\in L^2_\cF(t,T;\mathbb R^n)\times (L^2_\cF(t,T;\R^n))^d$
as  the unique solution to
\begin{equation} \label{adjoint1general}
\left\{\begin{array}{rl}
dp(s;t)= & \!\!\! -\bigg[A_s'p(s;t)+\displaystyle\sum_{j=1}^d(C_s^j)' k^j(s;t)+Q_sX^*_s\bigg]ds \\
             & +\displaystyle\sum_{j=1}^dk^j(s;t)dW_s^j,    \quad\quad s\in[t,T], \\
p(T;t)= & \!\!\! G X^*_T- h \Et{X^*_T}-\mu_1 X_t^*-\mu_2.
\end{array}\right.
\end{equation}

Notice that if $u^* \in L^p_\cF(0, \, T; \, K)$, then
$p(\cdot;t)\in L^p_\cF(t,T;\mathbb R^n)$ in fact.

Furthermore, define  $(P(\cdot;t), (K^j(\cdot;t))_{j=1,\cdots,d})\in L^\infty_\cF(t,T;\mathbb S^n)\times (L^2_\cF(t,T;\mathbb S^n))^d$ as  the unique solution to
\begin{equation}\label{adjoint2general}
\left\{\begin{array}{ll}
dP(s;t)=& \!\!\! -\bigg\{A_s'P(s;t)+P(s;t)A_s \\
&+\displaystyle\sum_{j=1}^d[(C_s^j)'P(s;t)C_s^j+(C_s^j)'K^j(s;t)+K^j(s;t)C_s^j]+Q_s\bigg\}ds \\
&+\displaystyle\sum_{j=1}^dK^j(s;t)dW_s^j, \quad\quad s\in[t,T], \\
P(T;t)=& \!\!\! G.
\end{array}\right.
\end{equation}

Notice that neither the terminal condition nor the coefficients of this equation depend on
$t$; so it can be taken as a BSDE on the entire time interval $[0,T]$. Denote its solution as
$(P(s), K(s))$, $s\in[0,T]$. It then follows from the
uniqueness of the solution to BSDE that $(P(s;t), K(s;t))=(P(s), K(s))$
at $s\in[t,T]$ for any $t\in[0,T]$.

The following estimate under local spike variation is reproduced from \cite[Proposition 3.1]{HJZ}.
\begin{proposition}\label{variate}
 For any $t\in [0,T)$, $\varepsilon>0$ and $v\in \displaystyle\bigcup_{p>2} L^p_\cF(0, \, T; \, K)$,  define $u^{t,\varepsilon,v}$
by (\ref{svgeneral}). Then
\begin{equation}\label{epsilongeneral}
J(t,X^*_t; u^{t,\varepsilon,v})-J(t,X^*_t; u^*)
=\mathbb E_t\int_t^{t+\varepsilon} \bigg(\ang{\Lambda(s;t),v_s-u^*_s}
+\frac{1}{2}\ang{H(s)(v_s-u^*_s), v_s-u^*_s}\bigg)ds+o(\varepsilon),
\end{equation}
where $\Lambda(s;t)\trieq   B_sp(s;t)+\sum_{j=1}^d(D_s^j)'k^j(s;t)+R_su^*_s$ and
$H(s)\trieq R_s+\sum_{j=1}^d (D_s^j)'P(s)D_s^j$.
\end{proposition}


In view of Proposition \ref{variate} and the fact that $H(s)\succeq 0$, it is straightforward to get the following sufficient condition of an equilibrium.
\begin{corollary}\label{equiv-lq}
A control $u^{*}\in L^{p}_\cF(0,T;K)$ is an equilibrium if
\begin{equation}\label{cond1}
\liminf_{\varepsilon\downarrow 0}\frac{1}{\varepsilon}\int_{t}^{t+\varepsilon}\Et{\left<\Lambda(s;t),v_s-u_s^*\right>}ds\ge 0,  \;\;a.s., \;\; \forall t\in [0,T).
\end{equation}
\end{corollary}

The necessary condition is somewhat different.
\begin{proposition}\label{necessary}
If a control $u^{*}\in L^{p}_\cF(0,T;K)$ is an equilibrium, then for $\theta\in (0,1]$,
\begin{equation}\label{necessary1}
\liminf_{\varepsilon\downarrow 0}\frac{1}{\varepsilon}\mathbb E_t\int_t^{t+\varepsilon} \left(\ang{\Lambda(s;t),v_s-u_s^*}
+\frac{\theta}{2}\ang{H(s) (v_s-u_s^*), (v_s-u_s^*)}\right)ds\ge 0.
\end{equation}
\end{proposition}
\pf
We set, for $v\in \displaystyle\bigcup_{p>2} L^p_\cF(0, \, T; \, K) $ and $\theta\in (0,1]$,
$$\bar{v}_s=u_s^*+\theta (v_s-u_s^*)\in L^{p}_\cF(t,T; K) .$$
Then
\begin{equation}\label{epsilontheta}
\begin{array}{rl}
& J(t,X^*_t; u^{t,\varepsilon,\bar{v}})-J(t,X^*_t; u^*) \\
=\!\!\! & \displaystyle\mathbb E_t\int_t^{t+\varepsilon} \!\!\! \bigg(\!\theta\ang{\Lambda(s;t),v_s-u_s^*}
+\frac{\theta^2}{2}\ang{H(s) (v_s-u_s^*), (v_s-u_s^*)}\!\bigg)ds+o(\varepsilon).
\end{array}
\end{equation}
Hence,
$$
\liminf_{\varepsilon\downarrow 0}\frac{1}{\varepsilon}\mathbb E_t\int_t^{t+\varepsilon} \left(\ang{\Lambda(s;t),v_s-u_s^*}
+\frac{\theta}{2}\ang{H(s) (v_s-u_s^*), (v_s-u_s^*)}\right)ds\ge 0.
$$
\eof

The next result provides a key property for the solution to
BSDE (\ref{adjoint1general}), and represents the process
$\Lambda(s;t)$ in a special form.

\begin{proposition}\label{k-ind-t}
For any given pair of state and control processes $(X^{*}, u^{*})$, the solution to BSDE (\ref{adjoint1general}) satisfies
$k(s;t_{1})=k(s;t_{2})$ for a.e. $s\ge \max{( t_{1},t_{2})}$.
Moreover, there exist $\lambda_{1}\in L^{p}_{\cF}(0, T; \R^{l}), \lambda_{2}\in L^{\infty}_{\cF}(0, T; \R^{l\times n}) $ and $\xi\in L^{p}(\Omega; C(0, T; \R^{n}))$, such that
$\Lambda(s;t)$ has the representation
$$\Lambda(s;t)=\lambda_{1}(s)+\lambda_{2}(s)\xi_{t}.$$
\end{proposition}
\pf
Define the function $\psi(\cdot)$ as the unique continuous solution to the following matrix-valued ordinary differential equation (ODE)
$$d\psi(t)=\psi(t)A(t)'dt, \;\;\;\; \psi(T)=I_{n},$$
where $I_{n}$ denotes the $n\times n$ identity matrix.  It is clear that
$\psi(\cdot)$ is invertible, and both $\psi(\cdot)$ and $\psi(\cdot)^{-1}$ are  bounded.

Let $\hat p(s;t)=\psi(s) p(s;t)+h\Et{X^{*}_{T}}+\mu_{1}X^{*}_{t}+\mu_{2}$
 and $\hat k^{j}(s;t)=\psi(s) k^{j}(s;t)$ for $j=1,\cdots, d$.
Then by It\^o's formula, on the time interval $[t, T]$,  $(\hat p(\cdot;t), \hat k(\cdot;t))$ satisfies
\begin{equation}\label{adjoint1tilde}
\left\{\begin{array}{l}
d\hat p(s;t)=-\left[\displaystyle\sum_{j=1}^d\psi(s)(C_s^j)' \psi(s)^{-1}\hat k^j(s;t)+\psi(s)Q_s X^*_s\right]ds+\displaystyle\sum_{j=1}^d\hat k^j(s;t)dW_s^j, \\
\hat p(T;t)=G X^*_T.
\end{array}\right.
\end{equation}

It is clear that neither the terminal condition nor the coefficients of this equation depend on
$t$; so it can be taken as a BSDE on the entire time interval $[0,T]$. We denote its solution as
$(\hat p(s), \hat k(s))$, $s\in[0,T]$. It  follows from the
uniqueness of the solution to BSDE that $(\hat p(s;t), \hat k(s;t))=(\hat p(s), \hat k(s))$
at $s\in[t,T]$ for any $t\in[0,T]$.
As a result, $k(s;t)=\psi(s)^{-1} \hat k(s):=k(s)$, proving the first claim of
the proposition.

Next, from the definition of $\hat p(s;t)$,
$$
p(s;t)=\psi(s)^{-1}\hat p(s)-\psi(s)^{-1}(h\Et{X^{*}_{T}}+\mu_{1}X^{*}_{t}+\mu_{2})=p(s)+\psi(s)^{-1}\xi_{t},
$$
 where 
 $\xi_{t}:=-h\Et{X^{*}_{T}}-\mu_{1}X^{*}_{t}-\mu_{2}$ defines the process $\xi\in L^{p}_{\cF}(\Omega; C(0, T; \R^{n}))$ and $p(s):=\psi(s)^{-1}\hat p(s)$ defines the process $p\in L^{p}_{\cF}(\Omega; C(0, T; \R^{n}))$. Hence,
 \begin{eqnarray*}
 \Lambda(s;t)&\!\!\!=\!\!\!&B_sp(s;t)+\sum_{j=1}^d(D_s^j)'k^j(s;t)+R_su^*_s\\
 &\!\!\!=\!\!\!&B_sp(s)+\sum_{j=1}^d(D_s^j)'k^j(s)+R_su^*_s+B_{s}\psi(s)^{-1}\xi_{t}\\
 &\!\!\!=\!\!\!&\lambda_{1}(s)+\lambda_{2}(s)\xi_{t},
\end{eqnarray*}
where $\lambda_{1}(s):=B_sp(s)+\sum_{j=1}^d(D_s^j)'k^j(s)+R_su^*_s$ and $\lambda_{2}(s):=B_{s}\psi(s)^{-1}$.
 \eof

We now set out to derive our general necessary and sufficient condition for equilibrium
controls. Although (\ref{cond1}) and (\ref{necessary1}) already provide  characterizing conditions, they are nevertheless not very useful
because they involve a limit. It is tempting to expect that the limit therein is $\Lambda(t;t)$,
in the spirit of the Lebesgue differentiation theorem. However,
one needs to be very careful, since in both  (\ref{cond1}) and (\ref{necessary1}), the conditional expectation with respect to ${\cal F}_t$ is
involved. 
The following lemma generalizes Lemma 3.4 in \cite{HJZ2} from $q=2$ to any $q>1$, and we provide
a complete proof here for the sake of self-containedness.

\begin{lemma}\label{cond-lebesgue}
Let $Y\in L^{q}_{\cF}(0,T; \R)$, $q>1$ be a given process. If $\displaystyle\liminf_{\varepsilon\downarrow 0}\frac{1}{\varepsilon}\int_{t}^{t+\varepsilon}\Et{Y_s}ds\ge 0,\;  a.e.  ~ t\in [0, T), a.s.$,
then $Y_t\ge 0,\;  a.e. ~ t\in [0, T), a.s.$.
\end{lemma}
\pf
Since $L^{q^*}_{\cF_{T}}(\Omega; \R_+)$ ($q^*$ is the conjugate of $q$) is a separable space, it follows from
the (deterministic) Lebesgue differentiation theorem that there is a countable dense subset  ${\mathcal D}\subset  L^{q^*}_{\cF_{T}}(\Omega;  \R_+)\cap L^{\infty}_{\cF_{T}}(\Omega;  \R_+)$,
such that for almost all $t$, we have
\begin{equation}\label{limtoprod}
\lim_{\varepsilon\downarrow 0}\frac{1}{\varepsilon}\int_{t}^{t+\varepsilon}\E{ Y_{s}\eta }ds=\E{ Y_{t}\eta }, \quad \forall \eta\in {\mathcal D},
\end{equation}
and $\displaystyle\lim_{\varepsilon\downarrow 0}\frac{1}{\varepsilon}\int_{t}^{t+\varepsilon}\E{Y_{s}^{q}}ds=\E{Y_{t}^{q}}.$

For any $\eta\in {\mathcal D}$, define $\eta_{s}={\mathbb E}_{s}[\eta]$. Then
$\E{ Y_{s}\eta }=\E{ Y_{s}\eta_s }$. We have the following estimates:
\begin{eqnarray*}
\left|\lim_{\varepsilon\downarrow 0}\frac{1}{\varepsilon}\int_{t}^{t+\varepsilon}\E{ Y_{s} (\eta_{s}-\eta_{t}) }ds\right|
&\!\!\!\le\!\!\!&\lim_{\varepsilon\downarrow 0}\frac{1}{\varepsilon}\left( \int_{t}^{t+\varepsilon}\E{Y_{s}^{q}}ds\right)^{\frac{1}{q}} \left(\int_{t}^{t+\varepsilon}\E{(\eta_{s}-\eta_{t})^{q^*}}ds \right)^{\frac{1}{q^*}}\\
&\!\!\!=\!\!\!&\lim_{\varepsilon\downarrow 0}\left(\frac{1}{\varepsilon} \int_{t}^{t+\varepsilon}\E{Y_{s}^{q}}ds\right)^{\frac{1}{q}}\left( \frac{1}{\varepsilon}\int_{t}^{t+\varepsilon}\E{(\eta_{s}-\eta_{t})^{q^*}}ds\right)^{\frac{1}{q^*}} \\
&\!\!\!\le\!\!\!&\lim_{\varepsilon\downarrow 0}\left(\frac{1}{\varepsilon} \int_{t}^{t+\varepsilon}\E{Y_{s}^{q}}ds \right)^{\frac{1}{q}}\left(\sup_{s\in [t,t+\varepsilon]}\E{(\eta_{s}-\eta_{t})^{q^*}}\right)^{\frac{1}{q^*}}\\
&\!\!\!\le\!\!\!&q\lim_{\varepsilon\downarrow 0}\left(\frac{1}{\varepsilon} \int_{t}^{t+\varepsilon}\E{Y_{s}^{q}}ds\right)^{\frac{1}{q}} \left(\E{(\eta_{t+\varepsilon}-\eta_{t})^{q^*}} \right)^{\frac{1}{q^*}}
=0,
\end{eqnarray*}
where the last inequality is due to Doob's martingale inequality as $\eta_s$ is an $L^{q^*}$-integrable martingale.
Hence for any $\eta\in {\mathcal D}$,
\begin{eqnarray*}
\E{ Y_{t}\eta_{t} }&\!\!\!=\!\!\!&\E{ Y_{t}\eta }\\
&\!\!\!=\!\!\!&\lim_{\varepsilon\downarrow 0}\frac{1}{\varepsilon}\int_{t}^{t+\varepsilon}\E{ Y_{s} \eta }ds\\
&\!\!\!=\!\!\!&\lim_{\varepsilon\downarrow 0}\frac{1}{\varepsilon}\int_{t}^{t+\varepsilon}\E{ Y_{s} \eta_{s} }ds\\
&\!\!\!=\!\!\!&\lim_{\varepsilon\downarrow 0}\frac{1}{\varepsilon}\int_{t}^{t+\varepsilon}\E{Y_{s}\eta_{t} }ds\\
&\!\!\!=\!\!\!&\lim_{\varepsilon\downarrow 0}\frac{1}{\varepsilon}\int_{t}^{t+\varepsilon}\E{(\Et{Y_{s}}) \eta_{t} }ds\\
&\!\!\!=\!\!\!&\lim_{\varepsilon\downarrow 0}\E{\left( \frac{1}{\varepsilon}\int_{t}^{t+\varepsilon}\Et{Y_{s}}ds\right) \eta_{t} }.
\end{eqnarray*}

Since (recall that $q^*$ is the conjugate of $q$)
\begin{eqnarray*}
\E{\left(\frac{1}{\varepsilon}\int_{t}^{t+\varepsilon}\Et{Y_{s}}ds\right)^{q}}
&\!\!\!=\!\!\!& \frac{1}{\varepsilon^q}\E{\left(\int_{t}^{t+\varepsilon}\Et{Y_{s}}ds\right)^{q}} \\
&\!\!\!\le\!\!\!& \frac{1}{\varepsilon^q}\E{\left(\int_{t}^{t+\varepsilon}ds\right)^{q/q^*}\int_{t}^{t+\varepsilon}\Et{Y_{s}}^{q}ds} \\
&\!\!\!=\!\!\!& \frac{1}{\varepsilon^q}\E{\varepsilon^{q/q^*}\int_{t}^{t+\varepsilon}\Et{Y_{s}}^{q}ds} \\
&\!\!\!=\!\!\!& \frac{1}{\varepsilon}\E{\int_{t}^{t+\varepsilon}\Et{Y_{s}}^{q}ds} \\
&\!\!\!\le\!\!\!&\frac{1}{\varepsilon}\int_{t}^{t+\varepsilon}\E{Y_{s}^{q}}ds,
\end{eqnarray*}
%
and $\displaystyle\lim_{\varepsilon\downarrow 0}\frac{1}{\varepsilon}\int_{t}^{t+\varepsilon}\E{Y_{s}^{q}}ds=\E{Y_t^q}$,
there exists a constant $\delta_t> 0$, such that
$$\E{\left(\frac{1}{\varepsilon}\int_{t}^{t+\varepsilon}\Et{Y_{s}}ds\right)^{q}}\le 2\E{Y_t^q}, \quad \forall\, \varepsilon\in (0, \delta_t).$$
This implies that $\frac{1}{\varepsilon}\int_{t}^{t+\varepsilon}\Et{Y_{s}}ds$ is uniformly integrable in
$\varepsilon\in (0, \delta_t)$.
Since $\eta$ is essentially bounded, so is $\eta_t$; hence by Fatou's lemma, for $a.e. ~ t\in [0,T]$, and  any $\eta \in {\mathcal D}$,
\begin{eqnarray*}
\E{ Y_{t}\eta }&\!\!\!=\!\!\!&\liminf_{\varepsilon\downarrow 0}\E{\left( \frac{1}{\varepsilon}\int_{t}^{t+\varepsilon}\Et{Y_{s}}ds\right) \eta_{t} }\\
&\!\!\!\ge\!\!\! &\E{\liminf_{\varepsilon\downarrow 0}\left(\left( \frac{1}{\varepsilon}\int_{t}^{t+\varepsilon}\Et{Y_{s}}ds\right) \eta_{t}\right) }\\
&\!\!\!\ge\!\!\!&0,
\end{eqnarray*}
which implies
$$Y_t\ge 0, \; a.e. ~ t\in [0,T], \; a.s..$$
\eof

\begin{theorem}\label{maingeneral}
Given a  control $u^*\in L^p_{\cF}(0,T; \R^l)$,
let $X^{*}$ be the corresponding state process
and $(p(\cdot; t),k(\cdot;t))\in L^p_{\cF}(t,T; \R^n)\times (L^2_{\cF}(t,T; \R^n))^{d}$
be the unique solution to BSDE (\ref{adjoint1general}).
Then $u^{*}$ is an equilibrium control if and only if
\begin{equation}\label{cond2}
\left<\Lambda(t;t),v_t-u^*_t\right>\ge 0, \mbox{a.s., a.e. } t\in [0, T].
\end{equation}
\end{theorem}
%
%
\pf
Recall that we have the representation $\Lambda(s; t)=\lambda_{1}(s)+\lambda_{2}(s)\xi_{t}$.
Then
\begin{eqnarray*}
& &\frac{1}{\varepsilon}\int_{t}^{t+\varepsilon}\Et{\left<\Lambda(s;t),v_s-u_s^*\right>}ds
-\frac{1}{\varepsilon}\int_{t}^{t+\varepsilon}\Et{\left<\Lambda(s;s),v_s-u_s^*\right>}ds\\
&\!\!\!=\!\!\!&\frac{1}{\varepsilon}\int_{t}^{t+\varepsilon}\Et{\left<\lambda_2(s)(\xi_t-\xi_s), v_s-u_s^*\right>}ds.
\end{eqnarray*}
Hence
\begin{equation}\label{difference}
\lim_{\varepsilon\downarrow 0}\left|\frac{1}{\varepsilon}\int_{t}^{t+\varepsilon}\Et{\left<\Lambda(s;t),v_s-u_s^*\right>}ds
-\frac{1}{\varepsilon}\int_{t}^{t+\varepsilon}\Et{\left<\Lambda(s;s),v_s-u_s^*\right>}ds\right|=0.
\end{equation}

If (\ref{cond2}) holds, then from (\ref{difference}),
$$\liminf_{\varepsilon\downarrow 0}\frac{1}{\varepsilon}\int_{t}^{t+\varepsilon}\Et{\left<\Lambda(s;t),v_s-u_s^*\right>}ds
=\liminf_{\varepsilon\downarrow 0}\frac{1}{\varepsilon}\int_{t}^{t+\varepsilon}\Et{\left<\Lambda(s;s),v_s-u_s^*\right>}ds\ge0,$$
i.e. (\ref{cond1}) holds, and from Corollary \ref{equiv-lq}, $u^*$ is an equilibrium.

Now we suppose that  $u^*$ is an equilibrium, then from (\ref{necessary1}) and (\ref{difference}),
$$\liminf_{\varepsilon\downarrow 0}\frac{1}{\varepsilon}\mathbb E_t\int_t^{t+\varepsilon} \left(\ang{\Lambda(s;s),v_s-u_s^*}
+\frac{\theta}{2}\ang{H(s) (v_s-u_s^*), (v_s-u_s^*)}\right)ds\ge 0.$$

Then, from Lemma \ref{cond-lebesgue}, for any $\theta\in (0,1]$,
$$\ang{\Lambda(t;t),v_t-u_t^*}
+\frac{\theta}{2}\ang{H(t) (v_t-u_t^*), (v_t-u_t^*)}\ge 0.$$
Sending $\theta\rightarrow 0^+$, we obtain (\ref{cond2}).
\eof

When $n=1$, the state process $X$ is a  scalar-valued process evolving by the dynamics
\begin{equation}\label{control}
dX_s=[A_sX_s+B_s'u_s+b_s]ds+ [C_sX_s+D_s u_s+\sigma_s]'dW_s;\quad X_0=x_0,
\end{equation}
where $A$ is a bounded deterministic  scalar function on $[0, T]$.
The other parameters $B,C,D$ are all essentially bounded  and $\cF_t$-adapted processes on $[0,T]$
with values in $\R^{ l}$,
$\R^d$, $\R^{ d\times l}$,  respectively. Moreover,  $b\in L^\infty_\cF(0,T; \R)$ and $\sigma\in L^\infty_\cF(0,T; \R^d)$.

In this case,  the two adjoint equations for the equilibrium become
\begin{eqnarray}
&&\left\{\begin{array}{l}
dp(s;t)=-[A_sp(s;t)+C'_s k(s;t)+Q_sX^*_s]ds+k(s;t)'dW_s,\;\; s\in[t,T], \\
p(T;t)=G X^*_T- h \mathbb E_t[X^*_T]-\mu_1 X_t^*-\mu_2;
\end{array}\right.  \label{adjoint1} \\
&&\left\{\begin{array}{l}
dP(s;t)=-[(2A_s+|C_s|^2) P(s;t)+2C'_sK(s;t)+Q_s]ds \\
\quad\quad\quad\quad\quad+\,K(s;t)'dW_s,\;\; s\in[t,T], \\
P(T;t)=G.
\end{array}\right. \label{adjoint2}
\end{eqnarray}

For  reader's convenience,  we state here the $n=1$ version of Theorem \ref{maingeneral}.
\begin{theorem}\label{main}
An admissible control  $u^*\in L^p_\cF(0,T;\R^l)$  is an equilibrium control if and only if,  for any time $t\in [0,T)$,
\begin{itemize}
 \item [{\rm (i)}] the  system of  SDEs
\begin{equation}\label{fbsde}
\left\{
\begin{array}{l}
dX^*_s=[A_sX^*_s+B'_su^*_s+b_s]ds+ [C_sX_s^*+D_su^*_s+\sigma_s]'dW_s,\;\; s\in[0,T], \\
X_0^*=x_0,\\
dp(s;t)=-[A_sp(s;t)+C'_sk(s;t)+Q_sX^*_s]ds+k(s;t)'dW_s,\;\; s\in [t,T], \\
p(T;t)=G X^*_T-h E_t[X^*_T]-\mu_1 X_t^*-\mu_2,\;\;t\in[0,T],
\end{array}\right.
\end{equation}
admits a solution $( X^*,p,k)$;
\item [{\rm (ii)}]  $\Lambda(\cdot; t)\trieq p(\cdot;t)B_{\cdot}+D'_{\cdot}k(\cdot;t)+R_{\cdot} u^*_{\cdot}$ satisfies the condition {\rm (\ref{cond2})}.
\end{itemize}
\end{theorem}

\section{Mean-Variance Equilibrium Strategies in a  Market under Convex Cone Constraint}\label{mvcase}

As an application of the time-inconsistent LQ theory, we consider the continuous-time Markowitz mean--variance portfolio selection
model in a  market under convex cone constraint with random model coefficients. We aim to establish the existence and uniqueness of the equilibrium strategy.
The model is mathematically a special case of the general LQ problem formulated
earlier in this paper, with $n = 1$ naturally.

We use the classical setup. The wealth equation is governed by the SDE
\begin{equation}\label{wealth2}
\left\{\begin{array}{l}
dX_s=[r_s X_s+\theta_s'u_s]ds+u_s'dW_s, \quad s\in [t,T],\\
X_t=x_t,
\end{array}\right.
\end{equation}
where $r$ is the (bounded) deterministic interest rate function, and $\theta$ is the essentially bounded stochastic
 risk premium process. In particular, $x_0>0$.

The objective at time $t$ with state $X_{t}=x_{t}$ is to minimize
\begin{eqnarray}\label{mv-obj}
J(t, x_t; u)&\!\!\!\trieq\!\!\! &\frac{1}{2}{\rm Var}_t(X_T)-\gamma_t(x_t) \mathbb E_t[X_T] \\
&\!\!\!=\!\!\!&\frac{1}{2}\left(\mathbb E_t[X_T^2]-(\mathbb E_t[X_T])^2\right)-\mu_1 x_t \mathbb E_t[X_T]. \nonumber
\end{eqnarray}
There are two sources of time-inconsistency in this model, one from the variance term and the other from the state-dependent tradeoff between
the mean and the variance. We suppose that the portfolio constraint $K$ is a convex cone here.

The FBSDE (\ref{controlgeneral:0}) and (\ref{adjoint1general}) specializes  to
\begin{equation}\label{fbsderand}
\left\{
\begin{array}{l}
dX^*_s=[r_s X^*_s+\theta_s' u^*_s]ds+ (u^*_s)'dW_s,\quad X_0^*=x_0, \\
dp(s;t)=-r_s p(s;t)ds+k(s;t)'dW_s, \\
p(T;t)=X^*_T-\mathbb E_t[X^*_T]-\mu_1 X_t^*,
\end{array}\right.
\end{equation}
and the process $\Lambda(s;t)$ in condition (\ref{cond2}) is
$$\Lambda(s;t)=p(s;t)\theta_s+k(s;t).$$

We require that
\begin{equation}\label{convex}
\left<\Lambda(t;t), v_t-u^*_t\right>\ge 0.
\end{equation}

\subsection{Existence}

In this subsection, we construct a solution to (\ref{fbsderand}) and (\ref{convex}).

Let us first assume the following Ansatz:
$$p(s;t)=M_s  X^*_s-\mathbb E_t[M_s X^*_s]-\rho_s\mu_1 X^*_t,$$
with
$$dM_s=-F_{M,U}(s)ds+U_s'dW_s,$$
and
$$\rho_s=e^{\int_s^Tr_vdv}.$$
Applying It\^o's formula to $M_sX^*_s$, we get
$$d(M_sX^*_s)=\big[(-F_{M,U}(s)+r_sM_s)X_s^*+(\theta_s M_s+U_s)'u_s\big]ds+(X_s^*U_s+M_su_s)'dW_s,$$
and then
$$d\mathbb E_t[M_sX^*_s]=\mathbb E_t\big[(-F_{M,U}(s)+r_sM_s)X_s^*+(\theta_s M_s+U_s)'u_s\big]ds.$$
Hence
$$k(s)=X_s^*U_s+M_su_s,$$
and
$$p(t;t)=-\rho_t\mu_1 X^*_t.$$
Then (\ref{convex}) becomes
$$\left<-\rho_t\mu_1\theta_t X^*_t+X_t^*U_t+M_tu_t, v_t-u^*_t\right>\ge 0.$$
We will construct a solution with $X_t^*\ge 0$, thus
$$u_t^*=\mbox{\rm Proj}_K \left( M_t^{-1}(\rho_t\mu_1\theta_t-U_t)\right)X_t^*.$$
Denoting
$$\alpha_t=\mbox{\rm Proj}_K \left( M_t^{-1}(\rho_t\mu_1\theta_t-U_t)\right),$$
and coming back to (\ref{fbsderand}), we get:
$$F_{M,U}(s)=2r_sM_s+(\theta_s M_s+U_s)'\alpha_s.$$

To proceed, let us recall some facts about bounded-mean-oscillation (BMO) martingales; see Kazamaki \cite{Kazamaki}.
The process $Z\cdot W\trieq\int_0^\cdot Z_s'dW_s$ is a BMO martingale if and only if there exists a constant $C>0$ such that
$$\mathbb E\left[\int_\tau^T |Z_s|^2ds\Big| {\cal F}_\tau \right] \le C$$
for any stopping time $\tau\le T$. For every such $Z$, the stochastic exponential of $Z\cdot W$ denoted by ${\cal E}(Z\cdot W)$ is a positive martingale,
and for any $p>1$, there exists a constant $C_p>0$ such that
$$\E{\left(\int_\tau^T|Z_s|^2ds\right)^p\Big|\cF_\tau}\le C_p$$
for any stopping time $\tau\le T$.
Moreover, if $Z\cdot W$ and $V\cdot W$ are both BMO martingales, then under the probability measure $\Q$
defined by $\frac{d\Q}{d\p}={\cal E}_T(V\cdot W)$, $W^\Q_t\trieq W_t-\int_0^tV_sds$ is a standard Brownian motion, and $Z\cdot W^\Q$ is a BMO martingale.

\begin{lemma}\label{q-BSDE}
The following quadratic BSDE
\begin{equation}\label{eq:lemma4.1}
dM_s=-\big[2r_sM_s+(\theta_s M_s+U_s)'\mbox{\rm Proj}_K \big( M_s^{-1}(\rho_s\mu_1\theta_s-U_s)\big)\big]ds+U_s'dW_s, Œ\quad M_T=1,
\end{equation}
admits a   solution $(M,U)\in L_{\cal F}^\infty(0,T;\mathbb R)\times L_{\cal F}^2(0,T;\mathbb R^d)$ satisfying  $M\ge c$
for some  constant $c>0$. Moreover, $U\cdot W$ is a BMO martingale.
\end{lemma}
\pf We can prove the existence by a truncation argument together with a Girsanov transformation.

Let $c>0$ be a given number to be chosen later. Consider the following quadratic BSDE:
\begin{equation}\label{MUtrun}
\left\{\begin{array}{l}
dM_s=-\big[2r_sM_s+\big(\theta_s (M_s\vee c) +U_s\big)'\mbox{\rm Proj}_K \big( (M_s\vee c)^{-1}(\rho_s\mu_1\theta_s-U_s)\big)\big] ds+U_s' dW_s, \\
M_T=1.\end{array}\right.
\end{equation}

This BSDE is a standard quadratic BSDE. Hence there exists a solution $(M^{c},U^{c})\in L_{\cF}^\infty(0,T;\R)\times L_{\cF}^2(0,T;\R^d)$,
and $U^{c}\cdot W$ is a BMO martingale; see \cite{Koby} and \cite{Morlais}.

We can rewrite the above BSDE as, by noticing that $\theta_s'\mbox{\rm Proj}_K(\theta_s)=|\mbox{\rm Proj}_K(\theta_s)|^2$ because $K$ is a convex cone,
\begin{equation}\label{MUtrun2}
\left\{\begin{array}{rl}
dM^c_s=\!\!\!&-(2r_s M^c_s+\rho_s\mu_1|\mbox{\rm Proj}_K(\theta_s)|^2)ds+(U^c_s)' [dW_s-\beta_sds], \\
M^c_T=\!\!\!&1,\end{array}\right.
\end{equation}
where
$$
\beta_s=\mbox{\rm Proj}_K\big( (M^c_s\vee c)^{-1}(\rho_s\mu_1\theta_s-U^c_s)\big)
+\frac{\theta_s'\big(\mbox{\rm Proj}_K( \rho_s\mu_1\theta_s-U^c_s)-\mbox{\rm Proj}_K(\rho_s\mu_1\theta_s)\big)}{|U^c_s|^2}U_s^c{\bf 1}_{\{U^c_s\not=0\}}.
$$

It is easy to see that $|\beta|\le C(1+|U^c|)$, hence $\beta\cdot W$ is a BMO martingale.

As $\beta\cdot W$ is a BMO martingale, there exists a new probability
measure $\mathbb Q$ such that
$$W_t^{\mathbb Q} = W_t-\int_0^t \beta_sds$$
is a Brownian motion under $\mathbb Q$.

Hence,
$$M_s^c=\mathbb E_s^{\mathbb Q}\left[ e^{2\int_s^T r_tdt}+\int_s^T \rho_v\mu_1e^{2\int_s^v r_tdt}|\mbox{\rm Proj}_K(\theta_v)|^2 dv\right],$$
from which we deduce that there exists a constant $\eta>0$ independent of $c$ such that $M\ge \eta$. Taking  $c= \eta$, we obtain a solution.
\eof

Now we can state our main existence theorem:
\begin{theorem}
The following feedback
$$u_s^*=\mbox{\rm Proj}_K\big( M_s^{-1}(\rho_s\mu_1\theta_s-U_s)\big)X_s^*$$
is an  equilibrium strategy.
\end{theorem}
\pf
Set
$$
\begin{array}{rcl}
dX^*_s & \!\!\!=\!\!\! & \big[r_s X^*_s+\theta_s'\mbox{\rm Proj}_K\big(M_s^{-1}(\rho_s\mu_1\theta_s-U_s)\big)X_s^*\big]ds \\
& & + \big[\mbox{\rm Proj}_K\big(M_s^{-1}(\rho_s\mu_1\theta_s-U_s)\big)X_s^*\big]'dW_s,\quad X_0^*=x_0; \\
u_s^* & \!\!\!=\!\!\! & \mbox{\rm Proj}_K\big(M_s^{-1}(\rho_s\mu_1\theta_s-U_s)\big)X_s^*, \\
p(s;t) & \!\!\!=\!\!\! & M_s  X^*_s-\mathbb E_t[M_s X^*_s]-\rho_s\mu_1 X^*_t,
\end{array}
$$
and
$$k(s)=X_s^*U_s+M_su^*_s.$$

Then, $X^*_s>0$, and
$$dX^*_s=[r_s X^*_s+\theta_s' u_s^*]ds+ (u_s^*)'dW_s,\quad X_0^*=x_0.$$

Let us now prove that $u^*$ is in $\displaystyle\bigcup_{p>2} L^p_\cF(0, \, T; \, K)$.
Applying Ito's formula to $M_s(X_s^*)^2$, we obtain
 (recall that $\alpha_s=\mbox{\rm Proj}_K\big(M_s^{-1}(\rho_s\mu_1\theta_s-U_s)\big)$),
\begin{eqnarray*}
d\big(M_s(X_s^*)^2\big) & \!\!\!=\!\!\! & -(X_s^*)^2\big[2r_sM_s+(\theta_s M_s+U_s)'\alpha_s\big]ds+(X_s^*)^2U_s'dW_s \\
         & &+\big(2M_s[r_s+\theta_s' \alpha_s]+M_s|\alpha_s|^2\big)(X_s^*)^2ds+ 2M_s\alpha_s'(X_s^*)^2dW_s+2(X_s^*)^2U_s'\alpha_sds \\
         & \!\!\!=\!\!\! & \bigg(\bigg[\theta_s+\frac{U_s}{M_s}\bigg]'\alpha_s+|\alpha_s|^2\bigg) M_s(X_s^*)^2ds+M_s(X_s^*)^2\bigg(2\alpha_s+\frac{U_s}{M_s}\bigg)'dW_s.
\end{eqnarray*}
As $K$ is a convex cone,
$$
a'\mbox{\rm Proj}_K(a)=|\mbox{\rm Proj}_K(a)|^2,
$$
from which we deduce that
$$
d\big(M_s(X_s^*)^2\big)= \big(1+M_s^{-1}\rho_s\mu_1\big)\theta_s'\alpha_s M_s(X_s^*)^2ds+M_s(X_s^*)^2\bigg(2\alpha_s+\frac{U_s}{M_s}\bigg)'dW_s.
$$
Hence
$$
M_t(X_t^*)^2 = M_0x_0^2e^{\int_0^t(1+M_s^{-1}\rho_s\mu_1)\theta_s'\alpha_sds}{\cal E}\bigg(\Big(2\alpha+\frac{U}{M}\Big)\cdot W\bigg)_t.$$

From John-Nirenberg's inequality (see Kazamaki \cite[Theorem 2.2, p.29]{Kazamaki}), we deduce that there exists $\varepsilon>0$ such that $\E{e^{\varepsilon \int_{0}^{T}|\alpha_{s}|^{2}ds}}<+\infty$.
Thus, $e^{\int_0^T(1+M_s^{-1}\rho_s\mu_1)\theta_s'\alpha_sds}\in \displaystyle\bigcap_{p>1} L^p$.

Moreover, as $(2\alpha+\frac{U}{M})\cdot W$ is a BMO martingale, $\displaystyle\sup_t\left[{\cal E} \bigg(\Big(2\alpha+\frac{U}{M}\Big)\cdot W\bigg)_t\right]\in \displaystyle\bigcup_{p>1} L^p$.

As $M\ge c>0$, we deduce that
$\displaystyle\sup_t|X_t^*|$ is in $\displaystyle\bigcup_{p>2} L^p$, and then $u^*$ is in $\displaystyle\bigcup_{p>2} L^p_\cF(0, \, T; \, K)$.

Now we calculate $dp$. Applying Ito's formula to $M_sX_s^*$, we obtain
\begin{eqnarray*}
d(M_sX_s^*)&\!\!\!=\!\!\!&-X_s^*\big[2r_sM_s+(\theta_s M_s+U_s)'\alpha_s\big]ds+X_s^*U_s'dW_s \\
         & &+M_s[r_s+\theta_s'\alpha_s]X_s^*ds+ M_s\alpha_s'X_s^*dW_s+X_s^*U_s'\alpha_sds \\
         &\!\!\!=\!\!\!&-r_sM_sX_s^*ds+k(s)'dW_s,
\end{eqnarray*}
and then
\begin{eqnarray*}
dp(s;t)&\!\!\!=\!\!\!&-r_sM_sX_s^*ds+k(s)'dW_s+r_s\mathbb E_t[M_sX_s^*]ds+r_s\rho_s\mu_1 X_t^*ds \\
       &\!\!\!=\!\!\!&-r_sp(s;t)ds+k(s)'dW_s.
\end{eqnarray*}
Hence $(X^*,u^*,p,k)$ is a solution to (\ref{fbsderand}), and (\ref{convex}) is easily checked.

\eof

\subsection{Uniqueness}

\begin{theorem} The following feedback
$$u_s^*=\mbox{\rm Proj}_K\big(M_s^{-1}(\rho_s\mu_1\theta_s-U_s)\big)X_s^*$$
is the unique  equilibrium strategy.
\end{theorem}

\pf
Suppose that $(X,u,p,k)$ is a solution to
\begin{equation}\label{fbsderand2}
\left\{
\begin{array}{l}
dX_s=[r_s X_s+\theta_s' u_s]ds+ (u_s)'dW_s,\quad X_0=x_0,\\
dp(s;t)=-r_s p(s;t)ds+k(s;t)'dW_s,\\
p(T;t)=X_T-\mathbb E_t[X_T]-\mu_1 X_t,
\end{array}\right.
\end{equation}
and the process $\Lambda(s;t)$ in condition (\ref{cond2}) is
$$\Lambda(s;t)=p(s;t)\theta_s+k(s;t).$$

We require that
\begin{equation}\label{convex2}
\left<\Lambda(t;t), v_t-u_t\right>\ge 0.
\end{equation}

There exist two adapted processes $\alpha_+$ and $\alpha_+$ with
$0\le \alpha_+\le 1$, $-1\le \alpha_-\le 0$ and $\alpha_+-\alpha_-=1$, such that
 $X^+_t=\alpha_+(t) X_t$, and $X^-_t=\alpha_-(t)X_t$.
Consider the following quadratic BSDE:
\begin{equation}\label{bsdetilde}
\begin{array}{rl}
d\tilde{M}_s=\!\!\!&-\bigg\{2r_s\tilde{M}_s+(\theta_s \tilde{M}_s+\tilde{U}_s)'\tilde{M}_s^{-1}\Big[\mbox{\rm Proj}_K \left((\rho_s\mu_1\theta_s-\tilde{U}_s)\right)\alpha_+(s) \\
&+\mbox{\rm Proj}_K \left(-(\rho_s\mu_1\theta_s-\tilde{U}_s)\right)\alpha_-(s)\Big]\bigg\}ds+\tilde{U}_s'dW_s.
\end{array}
\end{equation}

We note that
\begin{eqnarray*}
& &(\theta_s \tilde{M}_s)'\tilde{M}_s^{-1}\big[\mbox{\rm Proj}_K \left(\rho_s\mu_1\theta_s\right)\alpha_+(s)+\mbox{\rm Proj}_K \left(-\rho_s\mu_1\theta_s\right)\alpha_-(s)\big] \\
&=\!\!\!&\theta_s'\big[\mbox{\rm Proj}_K \left(\rho_s\mu_1\theta_s\right)\alpha_+(s)+\mbox{\rm Proj}_K \left(-\rho_s\mu_1\theta_s\right)\alpha_-(s)\big] \\
&=\!\!\!&\rho_s\mu_1[|\mbox{\rm Proj}_K(\theta_s)|^2\alpha_+(s)-|\mbox{\rm Proj}_K(-\theta_s)|^2\alpha_-(s)] \\
&\ge\!\!\!&0.
\end{eqnarray*}
Applying the same method as that of Lemma \ref{q-BSDE}, the quadratic BSDE (\ref{bsdetilde}) admits a  solution $(\tilde{M},\tilde{U})\in L_{\cal F}^\infty(0,T;\mathbb R)\times L_{\cal F}^2(0,T;\mathbb R^d)$ satisfying  $\tilde{M}\ge c$
for some  constant $c>0$. Moreover, $\tilde{U}\cdot W$ is a BMO martingale.
Let us take any such solution $(\tilde{M},\tilde{U})$.

It is important to note that
$$
\begin{array}{rl}
&\tilde{M}_{s}^{-1} \mbox{\rm Proj}_K \left((\rho_s\mu_1\theta_s-\tilde{U}_s)X_s\right) \\
=\!\!\! &\tilde{M}_{s}^{-1} \left[\mbox{\rm Proj}_K \left((\rho_s\mu_1\theta_s-\tilde{U}_s)\right)\alpha_+(s)+\mbox{\rm Proj}_K \left(-(\rho_s\mu_1\theta_s-\tilde{U}_s)\right)\alpha_-(s)\right]X_s \\
=\!\!\! & \tilde\alpha_sX_s,
\end{array}
$$
where
$$\tilde\alpha_s=\tilde{M}_{s}^{-1}[\mbox{\rm Proj}_K \left((\rho_s\mu_1\theta_s-\tilde{U}_s)\right)\alpha_+(s)+\mbox{\rm Proj}_K \left(-(\rho_s\mu_1\theta_s-\tilde{U}_s)\right)\alpha_-(s)].$$

Setting
$$
\bar{p}(s;t)=p(s;t)-\big(\tilde{M}_s{X}_s-\mathbb E_t[\tilde{M}_s{X}_s]-\rho_s\mu_1{X}_t\big),
$$
and
$$
\bar{k}(s)=k(s)-(X_s\tilde{U}_s+\tilde{M}_su_s).
$$
Then
$$\Lambda(s;t)=\big[\bar{p}(s;t)+\tilde{M}_sX_s-\mathbb E_t[\tilde{M}_sX_s]-\rho_s\mu_1{X}_t\big]\theta_s+\bar{k}(s)+X_s\tilde{U}_s+\tilde{M}_su_s,$$
and then condition (\ref{convex2}) becomes: for any $v_t\in K$,
$$\left<     [\bar{p}(t;t)-\rho_t\mu_1{X}_t]\theta_t+\bar{k}(t)+X_t\tilde{U}_t+\tilde{M}_tu_t, v_t-u_t \right>\ge 0,$$
from which we deduce that there exists one bounded adapted process $A$, such that
\begin{eqnarray*}
u_t&\!\!\!=\!\!\!&\tilde{M}_t^{-1}\mbox{\rm Proj}_K \left(-\bar{p}(t;t)\theta_t-\bar{k}(t)+(\rho_t\mu_1\theta_t-\tilde{U}_t)X_t)\right) \\
&\!\!\!=\!\!\!&\tilde{M}_t^{-1}A_t\left(-\bar{p}(t;t)\theta_t-\bar{k}(t)\right)+\tilde{M}_t^{-1}\mbox{\rm Proj}_K \left((\rho_t\mu_1\theta_t-\tilde{U}_t)X_t\right) \\
&\!\!\!=\!\!\!&\tilde{M}_t^{-1}A_t\left(-\bar{p}(t;t)\theta_t-\bar{k}(t)\right)+\tilde{\alpha}_tX_t.
\end{eqnarray*}

After some calculus, we arrive at:
\begin{equation}\label{barp0}
\left\{\begin{array}{rl}
d\bar p(s;t)=\!\!\!&-\Big\{ r_{s}\bar p(s;t)+(\theta_{s}+\tilde{U}_{s}\tilde{M}_{s}^{-1})'A_s[-\theta_{s}\bar p(s;s)-\bar k(s)] \\
& -\Et{(\theta_{s}+\tilde{U}_{s}\tilde{M}_{s}^{-1})'A_s[-\theta_{s}\bar p(s;s)-\bar k(s)]}\Big\}ds \\
& +\bar k(s)'dW_{s},\;\;s\in [t,T],\\
\bar p(T;t)=\!\!\!&0.
\end{array}\right.
\end{equation}

Applying the same method as in \cite{HJZ2}, we deduce that
$\bar p(s;t)=0$ and $\bar k(s)=0$. Therefore
$$u_t=\tilde{\alpha}_tX_t,$$
and then
$$X_t> 0, \quad \alpha_+=1, \quad \alpha_-=0.$$

As $X_t> 0$, replacing $(\tilde{M},\tilde{U})$ by $(M,U)$, and using the above procedure again,
we deduce then
$$u_s=\mbox{\rm Proj}_K\big(M_s^{-1}(\rho_s\mu_1\theta_s-U_s)\big)X_s,$$
and we conclude the proof.
\eof

From the above theorem, we deduce immediately the uniqueness of solution to BSDE (\ref{eq:lemma4.1}).

\begin{corollary}
The solution to BSDE (\ref{eq:lemma4.1}) is unique.
\end{corollary}

\pf
Let $(\bar M,\bar U)$ be another such solution.
Then from the above theorem, we deduce that
$$
\mbox{\rm Proj}_K\big(M_s^{-1}(\rho_s\mu_1\theta_s-U_s)\big)X^*_s=\mbox{\rm Proj}_K\big(\bar M_s^{-1}(\rho_s\mu_1\theta_s-\bar U_s)\big)X^*_s.
$$
As $X^*_s>0$, we deduce that,
$$
\mbox{\rm Proj}_K\big(M_s^{-1}(\rho_s\mu_1\theta_s-U_s)\big)=\mbox{\rm Proj}_K\big(\bar M_s^{-1}(\rho_s\mu_1\theta_s-\bar U_s)\big),
$$
from which we deduce the uniqueness of solution.
\eof

\subsection{Deterministic Risk Premium}
Let us first consider the case when the risk premium is a deterministic function of time.
Then $U=0$ and
$$M_s=e^{2\int_s^T r_vdv}\left( 1+\mu_1 \int_s^T e^{-\int_v^T r_zdz}|\mbox{\rm Proj}_K(\theta_v)|^2dv\right).$$
The equilibrium strategy is given by
$$u^*_s=\frac{\mu_1 e^{\int_s^T r_vdv}}{M_s}\mbox{\rm Proj}_K(\theta_s)X_s^*.$$

In the appendix, we obtain that the precommitted optimal control for the problem starting at $t=0$ is also in an affine feedback form
$$
u^{*pre}(s,x)=-\mbox{\rm Proj}_K(\theta_s)x+e^{\int_0^sr_vdv}\Big(x_0+\mu_1x_0 e^{\int_0^T(|\mbox{\footnotesize\rm Proj}_K(\theta_v)|^2-r_v)dv}\Big)\mbox{\rm Proj}_K(\theta_s).
$$

In \cite{BWYY}, the equilibrium is defined for the class of feedback controls as in \cite{BMZ}. Therein the equilibrium strategy is derived in a linear feedback form $u^{*fbe}_t=c_t^{fbe}X^*_t$ with $c^{fbe}_t$ uniquely determined by an integral equation and iterated by numerical method.
However, we can show that the linear coefficient of our equilibrium above is not a solution of the integral equation in \cite{BWYY}. This implies the difference between the two definitions of  equilibrium (open-loop and feedback).

To compare the performance of these two different equilibrium controls, together with the precommitted optimal control at time $t=0$,
we calculate $J(0,x_0;u)$ for $u=u^*, u^{*fbe}$, and $u^{*pre}$, respectively. Denote $c_s=\frac{\mu_1e^{\int_s^T r_vdv}}{M_s}\mbox{\rm Proj}_K(\theta_s)$,
then it is an easy exercise to get $X^*_T=x_0e^{\int_0^T(r_s+c_s'\theta_s-\frac{|c_s|^2}{2})ds+\int_0^Tc_s'dW_s}$. Hence
$$
\E{X^*_T}=x_0e^{\int_0^T(r_s+c_s'\mbox{\footnotesize\rm Proj}_K(\theta_s))ds}, \quad
\mbox{\rm Var}(X_T^*)=x_0^2e^{2\int_0^T(r_s+c_s'\mbox{\footnotesize\rm Proj}_K(\theta_s))ds}(e^{\int_0^T|c_s|^2ds}-1),
$$
leading to
$$
J(0,x_0;u^*)=\frac{x_0^2}{2}e^{2\int_0^T(r_s+c_s'\mbox{\footnotesize\rm Proj}_K(\theta_s))ds}\Big(e^{\int_0^T|c_s|^2ds}-1\Big) -x_0^2\mu_1 e^{\int_0^T(r_s+c_s'\mbox{\footnotesize\rm Proj}_K(\theta_s))ds}.
$$
Similarly,
$$
J(0,x_0;u^{*fbe})=\frac{x_0^2}{2}e^{2\int_0^T(r_s+(c_s^{fbe})'\mbox{\footnotesize\rm Proj}_K(\theta_s))ds}\Big(e^{\int_0^T|c^{fbe}_s|^2ds}-1\Big) -x_0^2\mu_1 e^{\int_0^T(r_s+(c_s^{fbe})'\mbox{\footnotesize\rm Proj}_K(\theta_s))ds}.
$$
By the calculation in the appendix, we have
$$
J(0,x_0;u^{*pre})=-\frac{x_0^2}{2}\mu_1^2\Big(e^{\int_0^T|\mbox{\footnotesize\rm Proj}_K(\theta_s)|^2ds}-1\Big)-x_0^2\mu_1e^{\int_0^Tr_sds}.
$$
Clearly,
$$
J(0,x_0;u^*)>J(0,x_0;u^{*pre}), \qquad  J(0,x_0;u^{*fbe})>J(0,x_0;u^{*pre}).
$$
Moreover, we can easily compare $J(0,x_0;u^*)$ and $J(0,x_0;u^{*fbe})$ due to their explicit expressions.

\section{Concluding Remarks}

In this paper, we consider some time-inconsistent LQ control problem under constraint.
We define the equilibrium strategy via spike perturbation of open control and deduce the necessary and
sufficient condition by applying the stochastic maximum principle, following the ideas of
\cite{HJZ, HJZ2}. LQ control problem with control constraint is useful because of its wide applications
in finance and economics. Our necessary and sufficient conditions are general enough to cover
many interesting time-inconsistent LQ control problem under various constraint.
We also shed light on important application in mean-variance portfolio under convex cone constraint and present its explicit equilibrium.
In particular, we can treat the random coefficient case, while the HJB method used by Bensoussan, Wong, Yam and Yung \cite{BWYY}
seems not applicable in random coefficient case.


\appendix
\section{Appendix. \\ Precommitted Mean--Variance Portfolio with Cone Constraint}

We consider the precommitted optimal control problem at time $t=0$,
\begin{equation}\label{premv2}
\begin{array}{cl}
\min & J(t, x_0; u) \trieq \displaystyle\frac{1}{2}{\rm Var}(X_T)-\gamma(x_0) \mathbb E[X_T] \\
&\quad\quad\quad\quad\quad \!\!=\displaystyle\frac{1}{2}\left(\mathbb E[X_T^2]-(\mathbb E[X_T])^2\right)-(\mu_1x_0+\mu_2)\mathbb E[X_T], \\
\mbox{s.t.} &dX_t=[r_tX_t+u_t'\theta_t]dt+u_t'dW_t,\\
&X_0=x_0.
\end{array}
\end{equation}

From the existing study on precommitted mean--variance problems such as \cite{HZ,LX,LZL}, it follows that, when the parameters $r_\cdot$ and $ \theta_\cdot$ are deterministic, we can get the explicit optimal value for the precommitted
problem (\ref{premv2}):
$$V^{pre}(x_0)=-\frac{1}{2}(\mu_1x_0+\mu_2)^2\Big(e^{\int_0^T|\mbox{\footnotesize\rm Proj}_K(\theta_s)|^2ds}-1\Big)-(\mu_1x_0+\mu_2)e^{\int_0^Tr_sds}x_0.$$
Furthermore, the corresponding optimal control can be written as the affine feedback control
$$u^{*pre}(s,x)=-\mbox{\rm Proj}_K(\theta_s)x+e^{\int_0^sr_vdv}\Big(x_0+(\mu_1x_0+\mu_2) e^{\int_0^T(|\mbox{\footnotesize\rm Proj}_K(\theta_v)|^2-r_v)dv}\Big)\mbox{\rm Proj}_K(\theta_s).$$

\end{document}